\documentclass{article}
\usepackage{amsmath,amsthm,amscd,amssymb,amsfonts}
\usepackage[all]{xy}
\UseComputerModernTips
\theoremstyle{plain}
	\newtheorem{theorem}{Theorem}[section]
	\newtheorem{lemma}[theorem]{Lemma}
	\newtheorem{proposition}[theorem]{Proposition}
	\newtheorem{corollary}[theorem]{Corollary}

\theoremstyle{definition}
	\newtheorem{example}[theorem]{Example}
	\newtheorem{definition}[theorem]{Definition}
	\newtheorem{remark}[theorem]{Remark}
\def\upabit{\vspace{-.3in} \\ \noindent}

\def\mapright#1{\buildrel #1 \over \longrightarrow}
\def\functor{\hspace{0.5mm}\underline{\hspace{2mm}}\hspace{0.5mm}}
\def\ab{{\rm ab}}

\def\CM{{\cal M}}

\def\CR{{\cal R}}

\def\C{\mathbb{C}}
\def\Z{\mathbb{Z}}
\def\N{\mathbb{N}}

\def\Coker{{\rm Coker}}
\def\Ker{{\rm Ker}}
\def\Im{{\rm Im}}

\def\Aut{{\rm Aut}}

\def\End{{\rm End}}

\def\GL{{\rm GL}}

\def\trace{{\rm Trace}}
\def\rEnd{\widetilde{\End}}
\def\rEndhzero{\widetilde{\End}^{\raisebox{-1.1mm}{\scriptsize{$h$}}}_0}
\def\ch{{\rm ch}}
\def\hch{\widehat{\ch}}
\title{
Non-commutative Characteristic Polynomials and Cohn Localization}
\author{Desmond Sheiham}
\begin{document}
\bibliographystyle{plain}
\maketitle
\begin{abstract}
Almkvist proved that for a commutative ring $A$ the characteristic
polynomial of an endomorphism $\alpha:P \to P$ of a finitely generated
projective
$A$-module determines $(P,\alpha)$ up to extensions.  For a
non-commutative ring $A$ the generalized characteristic polynomial of an
endomorphism $\alpha:P \to P$ of a finitely generated  projective $A$-module is
defined to be the Whitehead torsion $[1-x\alpha] \in K_1(A[[x]])$,
which is an equivalence class of formal power series with constant coefficient~$1$. 

In this paper an example is given of a non-commutative
ring $A$ and an endomorphism $\alpha:P \to P$ for which the
generalized characteristic polynomial does not determine $(P,\alpha)$
up to extensions.  The phenomenon is traced back to the non-injectivity
of the natural map $\Sigma^{-1}A[x] \to A[[x]]$,
where $\Sigma^{-1}A[x]$ is the Cohn localization of
$A[x]$ inverting the set $\Sigma$ of matrices in $A[x]$
sent to an invertible matrix by $A[x] \to A;x \mapsto 0$.
\end{abstract}
\section{Introduction}
We begin by recalling the definition of the characteristic
polynomial\footnote{ 
The polynomial defined here may be called the `reverse
characteristic polynomial' to distinguish between $\det(I-xM)$ and
$\det(M-xI)$. } 
$\ch_x(\C^n,\alpha)$ of an endomorphism $\alpha : \C^n\to \C^n$.
{\def\thefootnote{} \def\footnotesize{\normalsize} \footnotesep=5mm
\footnote{2000 Mathematics Subject Classification 16S34, 18F25.}}
\addtocounter{footnote}{-1}
\begin{equation*}
\ch_x(\C^n,\alpha)~=~\det(I - Mx)\in1+x\C[x]
\end{equation*}
where $M$ is an $n\times n$ matrix representing $\alpha$ with respect
to any choice of basis.

Of course, $\ch_x$ is not a complete invariant of the endomorphism;
for example the matrices
$\left(\begin{matrix}
\lambda & 0 \\
0 & \lambda
\end{matrix}\right)$ 
and~$\left(\begin{matrix}
\lambda & 1 \\
0 & \lambda
\end{matrix}\right)$
have the same characteristic polynomial although they are
not conjugate. On the other hand, if one is given the dimension~$n$ and the
characteristic polynomial $\ch_x(\C^n,\alpha)$, one can compute all the eigenvalues of~$\alpha$.
The Jordan normal form implies that
$(\C^n,\alpha)$ is determined uniquely up to choices of
extension (cf Kelley and Spanier~\cite{KelSpa68}).

The notion `unique up to choices of extension' can be made precise
without relying on a structure theorem for endomorphisms by
introducing the reduced
endomorphism class group $\rEnd_0(A)$ (see Almkvist~\cite{Alm73,Alm74} and Grayson~\cite{Gra77}) where $A$ denotes any ring.
$\rEnd_0(A)$ is\footnote{
Although free modules $A^n$ simplify the presentation,
the group $\rEnd_0(A)$ is unchanged if one substitutes
finitely generated projective modules throughout (see
section~\ref{subsection:Ktheorydefn}).
}
the abelian group with 
\begin{itemize}
\item one generator $[A^n,\alpha]$ for each isomorphism
class of pairs $(A^n,\alpha)$ where $\alpha:A^n\to A^n$, 
\item
a relation 
$[A^n,\alpha]+[A^{n''},\alpha''] = [A^{n'},\alpha']$
for each exact sequence
\begin{equation}
\label{equation:endoses}
0\to A^n \mapright{\theta} A^{n'} \mapright{\theta'} A^{n''} \to 0
\end{equation}
such that $\theta\alpha = \alpha'\theta$ and
$\theta'\alpha'=\alpha''\theta'$ and
\item
a relation $[A^n,0]=0$ for each $n$.
\end{itemize}
$\rEnd_0(\C)$, for example, is a free abelian group with one
generator $[\C,\lambda]$ for each non-zero eigenvalue
$\lambda\in\C\backslash 0$.

If $A$ is a commutative ring Almkvist proved~\cite{Alm74} that
the characteristic polynomial
\begin{equation}
\label{equation:commutative_ch}
\ch_x(A^n,\alpha)=\det(1-\alpha x:A[x]^n\to A[x]^n) 
\end{equation}
induces an isomorphism
\begin{equation*}
\ch_x:\rEnd_0(A) \to \widetilde{A}_0=\left\{\frac{1+a_1x+\cdots +a _nx^n}
{1+b_1x+\cdots+b_mx^m} \biggm{|} a_i,b_i\in A\right\}
\end{equation*}
so no further invariants are needed to classify endomorphisms up to
extensions.
 
If we do not assume that $A$ is commutative then the
definition~(\ref{equation:commutative_ch}) above does not
apply. However, $1-\alpha x:A[[x]]\to A[[x]]$ is a well-defined automorphism
(with inverse 
$1+\alpha x+\alpha^2x^2+\cdots$) where $A[[x]]$ denotes the ring of formal
power series in a central indeterminate $x$.  
 One can therefore define the {\it generalized characteristic polynomial} 
$\hch_x(A^n,\alpha)$ to be the element $[1-\alpha x]$ of the Whitehead group
$K_1(A[[x]])$, inducing a group homomorphism
\begin{equation*}
\hch_x : \rEnd_0(A) \to K_1(A[[x]]).
\end{equation*}
As Pajitnov observed~\cite{Paj95,PajRan00} a Gaussian
elimination argument (see section~\ref{subsection:pseries}) yields
\begin{equation*}
K_1(A[[x]])=K_1(A)\oplus W_1(A)
\end{equation*}
where $W_1(A)$ is the image in $K_1(A[[x]])$ of the group $1+xA[[x]]$
of Witt vectors.

If $A$ is commutative then $\widetilde{A}_0$ injects naturally 
into the group of units $A[[x]]^\bullet$ and
the commutative square 
\begin{equation*}
\xymatrix{\rEnd_0(A) \ar[r]^-{\displaystyle{\hch_x}} 
\ar[d]^-{\displaystyle{\cong}}_-{\displaystyle{\ch_x}}
& K_1(A[[x]]) 
\ar[d]^{\displaystyle{\det}} \\
\widetilde{A}_0 \ar[r] & A[[x]]^{\bullet}}
\end{equation*}
implies that $\hch_x$ is an injection.

The question arises whether $\hch_x$ is still injective when $A$ is
non-commutative. The main result of the present paper is that the answer can be
negative:   
\begin{proposition}
\label{Main_Proposition}
The non-commutative ring
\begin{equation*}
S=\Z\langle f,s,g\ |\ fg, fsg, fs^2g, \cdots\rangle.
\end{equation*}
is such that $\hch_x:\rEnd_0(S)\to K_1(S[[x]])$
is not injective.
\end{proposition}
Specifically the two endomorphisms $S\to S$
given by $a\mapsto as$ and $a\mapsto a(1-gf)s$ will be shown to have
the same image under $\hch_x$ although they represent distinct
classes in $\rEnd_0(S)$. The proof depends on the fact that the
functor $A\mapsto\rEnd_0(A)$ commutes with direct limits whereas 
$A\mapsto A[[x]]$ does not.

To put proposition~\ref{Main_Proposition} into context and explain the
origins of the ring $S$, we require a certain universal
localization $\Sigma^{-1}A[x]$ (Cohn~\cite[Ch.7]{Coh85},
Schofield~\cite[Ch4]{Sch85}) which P.M.Cohn constructed by adjoining
formal inverses to a set~$\Sigma$ of matrices.
Here, $\Sigma$ contains precisely the matrices which become invertible under 
the augmentation $\epsilon:A[x]\to A;x\mapsto 0$ (or equivalently are
invertible in $A[[x]]$). 

By the universal property of Cohn localization,
the inclusion of $A[x]$ in $A[[x]]$ factors in a unique way through
$\Sigma^{-1}A[x]$:  
\begin{equation}
\label{Universal_Triangle}
A[x]\mapright{i_\Sigma} \Sigma^{-1}A[x] \mapright{\gamma} A[[x]].
\end{equation}
In particular $i_\Sigma:A[x]\to \Sigma^{-1}A[x]$ is injective for
all rings $A$ (which is not true of some Cohn localizations).

If $A$ is commutative then $\Sigma^{-1}A[x]$ is the usual
commutative localization, inverting 
\begin{equation*}
\{\det(\sigma)\ |\ \sigma\in\Sigma\}=\{p\in A[x]\ |\ \epsilon(p)\
\text{is invertible}\},
\end{equation*} 
so $\gamma:\Sigma^{-1}A[x]\rightarrow A[[x]]$
is also injective. On the other hand
in section~\ref{section:gamma_not_injective} we prove:
\begin{proposition}
\label{Ring_Level_Prop}
The non-commutative ring $S$ is such that 
\begin{equation*}
\gamma:\Sigma^{-1}S[x] \to S[[x]]
\end{equation*}
is not injective.
\end{proposition}

In fact, proposition~\ref{Main_Proposition} is an algebraic $K$-theory version
of proposition~\ref{Ring_Level_Prop}; for a theorem due to Ranicki
(proposition~10.16 of~\cite{Ran98}) states that for any ring~$A$
\begin{equation}
\label{equation:Ranicki_Theorem}
K_1(\Sigma^{-1}A[x])\cong K_1(A)\oplus \rEnd_0(A)
\end{equation}
where the split injection $\rEnd_0(A)\to K_1(\Sigma^{-1}A[x])$ is 
$[A^n,\alpha]\mapsto [1-\alpha x]$. 
One can reinterpret proposition~\ref{Main_Proposition} as
the statement that the natural map $K_1(\Sigma^{-1}S[x]) \to
K_1(S[[x]])$ is not injective; by proposition~\ref{Ring_Level_Prop}
the phenomenon is not peculiar to algebraic $K$-theory. 

Propositions~\ref{Main_Proposition} and~\ref{Ring_Level_Prop} are
proved in sections~\ref{section:Ktheory}
and~\ref{section:gamma_not_injective} respectively. The proofs are
independent of each other and do not assume 
the identity~(\ref{equation:Ranicki_Theorem}) above.

Section~\ref{section:Many_Indeterminates} is expository. 
Firstly we 
show that, for any ring $A$, the image of $\gamma$ is the ring $\CR^A$
of rational power series; by definition $\CR^A$
is the smallest subring of $A[[x]]$ which contains
$A[x]$ and is such that elements of $\CR^A$ which are
invertible in $A[[x]]$ are invertible in $\CR^A$, i.e.~$\CR^A\cap
A[[x]]^\bullet = (\CR^A)^\bullet$. We work in greater 
generality replacing the
single indeterminate $x$ in~(\ref{Universal_Triangle}) by a set
$X=\{x_1,\cdots, x_\mu\}$ of non-commuting indeterminates
\begin{equation*}
A\langle X\rangle \mapright{i_\Sigma} \Sigma^{-1}A\langle X\rangle
\mapright{\gamma} A\langle\langle X\rangle\rangle.
\end{equation*}

Secondly we prove that each $\alpha\in\Sigma^{-1}A\langle
X\rangle$ can be expressed (non-uniquely) in the form
$\alpha=f(1-s_1x_1-\cdots-s_\mu x_\mu)^{-1}g$ where $f\in A^n$ is a
row vector, $g\in A^n$ is a column vector and 
$s_1,\cdots,s_\mu$ are $n\times n$ matrices with entries in~$A$. This
is a version of Sch\"utzenberger's 
theorem~\cite{Schu61, Schu62} (see also Berstel and
Reutenauer~\cite[Ch1]{BerReu88} and Cohn~\cite[\S6]{Coh75}). 
One can think of the elements of $\Sigma^{-1}A\langle X\rangle$ as
equivalence classes of finite dimensional linear machines
$(f,s_1,\cdots,s_\mu,g)$ which generate the power series
\begin{equation*}
%\label{equation:machine_output}
\gamma(\alpha)=fg+\sum_{i=1}^\mu fs_igx_i +
\sum_{i,j=1}^\mu fs_is_jgx_ix_j + \cdots~.
\end{equation*} 
Cohn wrote~\cite[p487]{Coh85}
\begin{quote}
The basic idea \ldots to invert matrices rather than elements was inspired by
the rationality criteria of Sch\"utzenberger and Nivat \ldots .
\end{quote}

Motivated by the theory of multi-dimensional boundary links, Farber
and Vogel proved~\cite{FarVog92} that if $A$ is a (commutative)
principal ideal domain then the Cohn localization of the free group
ring $AF_\mu$ (inverting those matrices which are invertible after
augmentation $AF_\mu\to A$) is isomorphic to the ring $\CR^A$ of
rational power series. In section~\ref{section:Farber_Vogel} we show
that this localization 
of the free group ring is isomorphic to $\Sigma^{-1}A\langle
X\rangle$ 
so $\gamma:\Sigma^{-1}A\langle X\rangle \to \CR^A$ is an isomorphism. 
By contrast, proposition~\ref{Ring_Level_Prop} above says that
$\Sigma^{-1}S\langle X\rangle$ is larger than $\CR^S$ even when $|X|=1$;
 distinct classes of linear machines can generate the same rational
power series. 
\bigskip  

I would like to thank my PhD supervisor
Professor Andrew Ranicki for all his help and encouragement.
%%%
%%%
\section{Algebraic $K$-theory}
\label{section:Ktheory}
\subsection{Definitions}
\label{subsection:Ktheorydefn}
Let $A$ be a ring, assumed to be associative and to contain a $1$.
We recall first the definitions of the Grothendieck group $K_0(A)$,
the Whitehead group $K_1(A)$ 
and the less widely known endomorphism class group 
\begin{equation*}
\End_0(A)=K_0(\text{Endomorphism category over A}).
\end{equation*}
\begin{definition}
$K_0(A)$ is the abelian group with one generator $[P]$ for each isomorphism
class of finitely generated projective $A$-modules and one
relation $[P']=[P]+[P'']$ for each identity $P'\cong P\oplus P''$.
\end{definition}

Let $\End(A)$ denote the category of pairs $(P,\alpha)$ where $P$ is a
projective (left) $A$-module and $\alpha:P\to P$ is an $A$-module
endomorphism. A morphism $\theta:(P,\alpha)\to (P',\alpha')$ in
$\End(A)$ is an $A$-module map $\theta:P\to P'$ such that $\theta \alpha =
\alpha'\theta$. A sequence of objects and morphisms
\begin{equation}
\label{equation:projendoses}
0\to (P,\alpha)\mapright{\theta} (P',\alpha') \mapright{\theta'}
(P'',\alpha'')\to 0
\end{equation}
is exact if $0\to P\mapright{\theta} P' \mapright{\theta'}
P''\to 0$ is an exact sequence.

Let $\Aut(A)\subset \End(A)$ denote the full subcategory of pairs 
$(P,\alpha)$ such that $\alpha:P\to P$ is an automorphism. 
\begin{definition}
The Whitehead group $K_1(A)$ is
the abelian group generated by the isomorphism classes $[P,\alpha]$ of
$\Aut(A)$ subject to relations: 
\begin{enumerate}
\item If $0\to (P,\alpha)\to (P',\alpha') \to
(P'',\alpha'')\to 0$ is an exact sequence  then
$[P',\alpha']=[P'',\alpha'']+[P,\alpha]$.
\item $[P,\alpha]+[P,\alpha']=[P,\alpha\alpha']$.
\end{enumerate}
\end{definition}
Alternatively, in terms of matrices, 
\begin{equation*}
K_1(A)=\GL(A)^{{\rm ab}}=\frac{\GL(A)}{E(A)}=\frac{\varinjlim
\GL_n(A)}{\varinjlim E_n(A)}
\end{equation*} 
where $E_n(A)$ is the subgroup of $\GL_n(A)$ generated by elementary
matrices $e_{ij}(a)$ which have $1$'s on the diagonal, $a$ in the
$ij$th position and $0$'s elsewhere ($a\in A$, $1\leq i,j\leq n$ and $i\neq
j$). See for example Rosenberg~\cite{Ros94} for further details. If
$M,M'\in \GL(A)$ and $[M]=[M']\in K_1(A)$ then we write $M\sim M'$. 

\begin{definition}
\label{definition:End_0(A)}
The endomorphism class group $\End_0(A)=K_0(\End(A))$ is the abelian group
with one generator $[P,\alpha]$ for each isomorphism class in $\End(A)$ and
a relation
\begin{equation}
\label{additive_relation}
[P',\alpha']=[P'',\alpha''] + [P,\alpha]
\end{equation}
corresponding to each exact sequence~(\ref{equation:projendoses}) above.
\end{definition}

Since every exact sequence of projective modules splits, 
we recover $K_0(A)$ by omitting the endomorphisms in
definition~\ref{definition:End_0(A)}. 
The forgetful map 
\begin{equation*}
\End_0(A)\to K_0(A);\ [P,\alpha]\mapsto [P]
\end{equation*}
is surjective and split by $[P]\mapsto [P,0]$ so that
$\End_0(A)\cong K_0(A)\oplus \rEnd_0(A)$
with
\begin{equation*}
\rEnd_0(A)~=~\Ker(\End_0(A)\to K_0(A))~\cong~\Coker(K_0(A)\to
\End_0(A))~.
\end{equation*}
Note that $\End_0(\functor)$ and $\rEnd_0(\functor)$ are functors; a
ring homomorphism $p:A\to A'$ induces a group homomorphism
\begin{align*}
\End_0(A) &\to \End_0(A') \\
[P,\alpha] &\mapsto [A'\otimes_A P, 1\otimes \alpha]. 
\end{align*}

Lemma~\ref{lemma:freemodules} below shows that the same group
$\rEnd_0(A)$ is obtained if, as in the introduction, one starts with
free modules in place of projective modules. Let 
$K_0^h(A)$ denote the Grothendieck group generated by free modules
$[A^n]$ subject to relations $[A^{m+n}]=[A^m]+[A^n]$.
Nearly all of the rings usually
encountered (including the ring~$S$ of the present paper) have
`invariant basis number', $A^n\cong A^m \Rightarrow n=m$, which implies
$K^h_0(A)=\Z$. 

Let $\End^h(A)\subset \End(A)$ denote the full
subcategory of pairs $(A^n,\alpha)$. Then $\End^h_0(A)=K_0(\End^h(A))$
satisfies 
$\End^h_0(A)\cong K_0^h(A)\oplus \rEndhzero(A)$ where
\begin{equation*}
\rEndhzero(A)~=~\Ker(\End^h_0(A)\to K_0^h(A))~\cong~\Coker(K_0^h(A)\to
\End_0^h(A))~.
\end{equation*}
\begin{lemma}
\label{lemma:freemodules}
There is a natural isomorphism $\rEndhzero(A)\cong\rEnd_0(A)$.
\end{lemma}
\begin{proof}
The homomorphism 
\begin{align*}
\rEndhzero(A) \cong \frac{\End^h_0(A)}{\Z} &\to
\frac{\End_0(A)}{K_0(A)} \cong\rEnd_0(A) \\
[A^n,\alpha] &\to [A^n,\alpha] 
\end{align*}
has inverse $[P,\alpha] \mapsto [P\oplus Q,\alpha\oplus 0]$
where $Q$ is a finitely generated $A$-module such that
$P\oplus Q$ is free. The definition of the inverse does not depend on
the choice of $Q$ and plainly $[P,0]\mapsto 0$ so we need only check
that the `exact sequence relations'~(\ref{additive_relation}) are respected. 
Suppose we are given an exact
sequence~(\ref{equation:projendoses}). Choose finitely generated $A$-modules $Q$ and $Q''$
 such that  $P\oplus Q$ and
$P''\oplus Q''$ are free. Then $P'\oplus Q\oplus
 Q''\cong P\oplus 
P'' \oplus Q\oplus Q''$ is free and there is an exact sequence of
endomorphisms:
\begin{equation*}
0\to (P\oplus Q,\alpha\oplus 0)\to (P'\oplus Q\oplus
Q'',\alpha'\oplus 0\oplus 0)\to (P''\oplus Q'',\alpha\oplus
0)\to 0~.
\end{equation*}
\upabit
\end{proof}

It follows from lemma~\ref{lemma:freemodules} that $\rEnd_0(A)$ has
an equivalent definition in terms of matrices:
let $M_n(A)$ denote the ring of $n\times n$ matrices with entries in $A$.
Regarding $A^n$ as a module of row vectors, a matrix $M\in M_n(A)$
represents the endomorphism of $A^n$ which multiplies by $M$ on
the right.
$\rEnd_0(A)$ is isomorphic to the group generated by $\left\{[M]\ \big{|}\
M\in\bigcup_{n=1}^\infty M_n(A)\right\}$ 
subject to relations:
\begin{enumerate}
\item
If $M\in M_n(A)$ and $M'\in M_{n'}(A)$ then
$[M]+[M']=\left[
\begin{matrix}
M & N \\
0 & M'
\end{matrix}
\right]$ for all $n\times {n'}$ matrices $N$.
\item
If $M,P\in M_n(A)$ and $P$ is invertible then $[M]=[PMP^{-1}]$.
\item
If all the entries in $M$ are zero then $[M]=0$.
\end{enumerate}

\subsection{Rings of Formal Power Series}
\label{subsection:pseries}
Let $A[[x]]$ be the ring of formal power series in the central
indeterminate $x$.

To define the generalized characteristic
polynomial of $(P,\alpha)$ we observe that $1-\alpha x$ has inverse
$1+\alpha x + \alpha^2 x^2 + \cdots$ when regarded as an
endomorphism of $P[[x]]=A[[x]]\otimes_A P$. Thus $1-\alpha x$ represents an
element of $K_1(A[[x]])$. 
Now
\begin{align*}
\hch_x : \rEnd_0(A) &\to  K_1(A[[x]]) \\
[P,\alpha] &\mapsto [1-\alpha x : P[[x]]\to P[[x]]\ ]
\end{align*}
is well defined because $\hch_x(P,0)=0\in K_1(A[[x]])$
and an exact sequence~(\ref{equation:projendoses}) gives rise to an exact
sequence
\begin{equation*}
0\to (P[[x]],1-\alpha x)\to (P'[[x]],1-\alpha'x)\to
(P''[[x]],1-\alpha''x)\to 0~.
\end{equation*}
\begin{lemma}
\label{lemma:K1powerseries}
i) $K_1(A[[x]])=K_1(A)\oplus W_1(A)$ 
where $W_1(A)$ is the image in $K_1(A[[x]])$ of
the group $W(A)=1+xA[[x]]$. \smallskip \\ \noindent
ii) If $A$ is commutative then $W_1(A)=W(A)=1+xA[[x]]$.
\end{lemma}
\noindent
This result and an argument 
showing that the abelianized group $(1+xA[[x]])^\ab$ is in general
larger than $W_1(A)$ can be found in Pajitnov and Ranicki~\cite{PajRan00}.
 
\begin{proof}[Proof of Lemma]
i) Let $\epsilon$ denote the augmentation map $A[[x]]\rightarrow A;x\mapsto0$. 
We shall prove that the sequence
\begin{equation*}
0\rightarrow W_1(A)\rightarrow K_1(A[[x]])\mapright{\epsilon}
K_1(A)\rightarrow 0 
\end{equation*}
is split exact. 

The composite $A \to A[[x]]\mapright{\epsilon} A$ is the
identity map so $\epsilon:K_1(A[[x]])\rightarrow K_1(A)$ is surjective and
split.

We have only to show that an element $\delta$ of
$K_1(A[[x]])$ which becomes zero in $K_1(A)$ can be written
$\delta=[1+x\xi]$ for some $\xi\in A[[x]]$. 
We may certainly write $\delta=[\delta_0 + x\delta_1 + x^2\delta_2
+\cdots]$ with 
$\delta_i\in M_n(A)$ for each $i$ and with $\delta_0$ invertible.
Now $[\delta_0]=0\in K_1(A)$ so
$\delta= [1+\eta ]$ where $\eta=\sum_{i=1}^\infty
\delta_0^{-1}\delta_ix^i$.
Since the diagonal entries of $1+\eta$ are invertible and all other entries are in
$xA[[x]]$, we can reduce $1+\eta$ by
elementary row operations to a diagonal matrix with entries
in $1+xA[[x]]$. Thus $\delta=[1+x\xi]$ where $1+x\xi$ is the product
of the diagonal entries. \smallskip \\ \noindent
ii) Taking determinants gives a homomorphism to the group of units
\begin{equation*}
\det:K_1(A[[x]])\rightarrow A[[x]]^\bullet
\end{equation*}
Every element of $W_1(A)$ can be written in the form $[1+x\xi]$ 
so the restriction of $\det$ to $W_1(A)$ is inverse to
the canonical map $1+xA[[x]]\rightarrow W_1(A)$.
\end{proof}
\subsection{Proof of Proposition~\ref{Main_Proposition}}
Recall that $S$ denotes the quotient of
the free ring $\Z\langle f,s,g\rangle$ by the two sided ideal
generated by the set $\{fs^ig\ |\ i=0,1,2,\cdots\}$.
There are two statements to prove:
\begin{itemize}
\item[\ref{Main_Proposition}a)]  $[S,s]$ and $[S, (1-gf)s]$ are
distinct classes in $\rEnd_0(S)$. 
\item[\ref{Main_Proposition}b)] $\hch_x[S,s]=\hch_x[S,(1-gf)s]$ in $K_1(S[[x]])$.
\end{itemize}
\noindent
\begin{proof}[Proof of {\rm \ref{Main_Proposition}b)}]
We aim to show $[1-sx]=[1-(1-gf)sx] \in K_1(S[[x]])$.
In $S[[x]]$ we have
\begin{equation*}
f(1-sx)^{-1}g~=~fg+(fsg)x+(fs^2g)x^2+(fs^3g)x^3 +\cdots~=~0~
\end{equation*}
so
\begin{align*}
1-sx &\sim
\left(\begin{matrix}
1 & 0 \\
0 & 1-sx
\end{matrix}\right)
\left(\begin{matrix}
1 & f \\
0 & 1
\end{matrix}\right)
\left(\begin{matrix}
1+f(1-sx)^{-1}g & 0 \\
0 & 1
\end{matrix}\right)
\left(\begin{matrix}
1 & 0 \\
-(1-sx)^{-1}g & 1
\end{matrix}\right) \\
&=
\left(\begin{matrix}
1 & f \\
-g & 1-sx
\end{matrix}\right) \\
&\sim
\left(\begin{matrix}
1 & f \\
0 & 1 
\end{matrix}\right)
\left(\begin{matrix}
1 & 0 \\
g & 1 
\end{matrix}\right)
\left(\begin{matrix}
1 & -f \\
0 & 1 
\end{matrix}\right)
\left(\begin{matrix}
1 & f \\
-g & 1-sx
\end{matrix}\right)
\left(\begin{matrix}
1 & -f \\
0 & 1
\end{matrix}\right) \\
&=
\left(\begin{matrix}
1 & 0 \\
0 & 1-(1-gf)sx
\end{matrix}\right) \ \ \ \mbox{since $fg=0$} \\
&\sim 1-(1-gf)sx.
\end{align*}
\end{proof}

To prove \ref{Main_Proposition}a), it is convenient to define a 
second invariant $\chi$. In terms of matrices,
\begin{align*}
\chi:  \rEnd_0(A) &\to \overline{A}[[x]] \\
[M] &\mapsto \sum_{i=1}^\infty \trace(M^i)x^i
\end{align*}
where $\overline A$ denotes the quotient of $A$ by the abelian group 
generated by commutators~(cf~Pajitnov~\cite{Paj00})
\begin{equation*}
\overline A = \frac{A}{\Z\{ab-ba\ |\ a,b\in A\}}.
\end{equation*}
\begin{example}
\label{example:free_ring_commutators}
Let $X$ be a set and suppose $A$ is the free ring
$\Z\langle X\rangle$ generated by $X$. The free monoid $X^*$ of words
in the alphabet $X$ is a basis for $\Z\langle X\rangle$ as a
$\Z$-module.
Each commutator $ab-ba$ with $a,b\in \Z\langle X\rangle$
is a linear 
combination of `basic' commutators $\sum_i\lambda_i (u_iv_i-v_iu_i)$
where $\lambda_i\in \Z$ and $u_i,v_i\in X^*$ so the
commutator submodule $\Z\{ab-ba\ |\ a,b\in\Z\langle X\rangle\}\subset
\Z\langle X\rangle$ is
spanned by elements $w-w'$ with  
$w,w'\in X^*$ and $w'$ a cyclic permutation of $w$ (written 
$w\sim w'$). Thus $\overline{\Z\langle X\rangle}=\Z\{X^*/\sim\}$.
\end{example} 
We emphasize that the abelian group $\overline A$ is in general 
larger than the commutative ring $A^{{\rm ab}}$, the latter being the
quotient of $A$ by the two-sided ideal generated by $\{ab-ba\ |\ a,b\in A\}$.
Nevertheless, if $M$ and $N$ are $n\times n$ matrices 
with
entries in $A$ then $\trace(MN)=\trace(NM)\in \overline A$ and it follows that $\chi$
is well-defined on $\rEnd_0(A)$.

\begin{remark}
$\chi$ is in general a weaker invariant then
$\hch_x$. There is a commutative triangle
\begin{equation*}
\xymatrix{
\rEnd_0(A)\ar[r]^-{\displaystyle{\hch_x}} \ar[rd]^-{\displaystyle{\chi}} & K_1(A[[x]]) \ar[d]^-{\displaystyle{T}} \\
& \overline A[[x]]
}
\end{equation*}
where 
\begin{equation*}
T[M]=-\trace\left(\left(x\frac{d}{dx}M\right)M^{-1}\right)
\end{equation*}
for $M\in \GL(A[[x]])$. 
Differentiation is defined formally by
\begin{equation*}
\frac{d}{dx}\sum_{n=0}^\infty a_nx^n = \sum_{n=1}^\infty na_nx^{n-1}~.
\end{equation*}
\end{remark}
\begin{proof}[Proof of {\rm \ref{Main_Proposition}a)}]
We define a family of rings
\begin{equation*}
S_m:=\Z\langle f,s,g\ |\ fg,fsg,\cdots,fs^mg\rangle\ .
\end{equation*}
There is an obvious surjection $p_m:S_m\twoheadrightarrow S_{m+1}$ for
each $m\in\N$ and $S$ is the direct limit of the system
%\begin{equation*}
$S=\varinjlim S_m.$
%\end{equation*}

By lemma~\ref{lemma:rEnd_directlimits} of appendix~\ref{appendix:directlimits}
we have $\rEnd_0(S)=\varinjlim \rEnd_0(S_m)$ so 
it suffices to prove that for each $m\in \N$
\begin{equation*}
[S_m,s]\neq [S_m,(1-gf)s]  \in \rEnd_0(S_m)~.
\end{equation*}

We shall see that $\chi$ is sensitive enough to distinguish these two
endomorphism classes. Indeed,
$\chi[S_m,(1-gf)s] = \sum_{i=1}^\infty ((1-gf)s)^ix^i$ and in particular
the coefficient of $x^{m+1}$ is   
\begin{equation*}
((1-gf)s)^{m+1}=s^{m+1} - (gfs^{m+1} + sgfs^{m} + \cdots + s^{m}gfs)
+ \text{other terms}
\end{equation*}
\vspace{-.3in} \\ \noindent
 where in each of the `other terms' two or more occurrences of $gf$
intersperse $m+1$ copies of $s$.
 Since $ab=ba \in \overline{S_m}$ for all $a,b\in S_m$, one may perform
a cyclic permutation of the letters in each term to obtain 
\begin{equation*}
((1-gf)s)^{m+1}=s^{m+1} - (m+1)fs^{m+1}g,
\end{equation*}
the `other terms' disappearing by the defining relations
$fg=\cdots=fs^mg=0$ of~$S_m$.
Now the coefficient of $x^{m+1}$ in $\chi[S,s]$ is $s^{m+1}$ so
it remains to prove that $(m+1)fs^{m+1}g \neq 0$ in
$\overline{S_m}$. We shall argue by contradiction.

Let $X$ denote the alphabet $\{f,s,g\}$. 
If $(m+1)fs^{m+1}g =0\in \overline{S_m}$ then there is an equation in
$\Z\langle X\rangle$:
\begin{equation}
\label{equation:absurd_commutation}
(m+1)fs^{m+1}g = \sum_{i=1}^l (w_i-w'_i) + r_0fg r'_0 + r_1 fsg r'_1 +
\cdots + r_m fs^m g r'_m
\end{equation}
where $r_j,r'_j \in \Z\langle X\rangle$ for $1\leq j\leq m$ and
$w_i,w'_i\in X^*$ are such that $w_i\sim w'_i$ for $1\leq i\leq l$ as
in example~\ref{example:free_ring_commutators} above. 

Let $V$ denote the $\Z$-module
generated by the cyclic permutations of $fs^{m+1}g$ and
let $W$ be the $\Z$-module generated by all other words in $X^*$
\begin{equation*}
\Z\langle X\rangle= V\oplus W=\Z\{w\in X^*\ |\ w\sim
fs^{m+1}g\}\oplus\Z\{w\in X^*\
|\ w\nsim fs^{m+1}g\}.
\end{equation*}
\vspace{-.3in} \\ \noindent
Each basic commutator $w-w'$ is either in $V$ or in $W$ and
\begin{equation*}
r_0fgr'_0+r_1fsgr'_1 + \cdots + r_mfs^mgr'_m \in W
\end{equation*}
so by equation~(\ref{equation:absurd_commutation})  
\begin{equation*}
(m+1)fs^{m+1}g=\sum_{i\in I}(w_i-w'_i)
\end{equation*}
 where $I=\{i\ |\ w_i \sim fs^{m+1}g\}\subset \{1,\cdots,l\}$.
We have reached a contradiction (for example put~$f=g=s=1$) and the proof of
proposition~\ref{Main_Proposition} is complete.
\end{proof}
%%%
%%%
\section{Cohn Localization}
\label{section:gamma_not_injective}
In this section, we briefly review Cohn localization
before proving proposition~\ref{Ring_Level_Prop}. 

\subsection{Definitions}
If $A$ is a ring and $\Sigma$ is any set of
matrices with entries in $A$ then a ring homomorphism $A\to B$ is said to
be $\Sigma$-inverting if every matrix in $\Sigma$ is mapped to an
invertible matrix over $B$. The Cohn localization 
$i_\Sigma:A\to \Sigma^{-1}A$ is the (unique) ring homomorphism with
the universal property that every $\Sigma$-inverting homomorphism
$A\to B$ factors uniquely through $i_\Sigma$. 
Note that $i_\Sigma$ is not in general an
injection; it may even be the case that $\Sigma^{-1}A=0$.

If $A$ is commutative then
$\Sigma^{-1}A$ coincides with the commutative ring of quotients $S^{-1}R$ with
$S=\{\det(M)\ |\ M\in\Sigma\}$.

For non-commutative $A$, Cohn constructed $\Sigma^{-1}A$ by generators
and relations as follows~\cite[p390]{Coh85}. For each $m\times n$
matrix $M\in\Sigma$ take a set of $mn$ symbols arranged as an $n\times
m$ matrix $M'$. $\Sigma^{-1}A$ is generated by the elements of $A$ together
with all the symbols in the matrices $M'$, subject to the relations
holding in $A$ and the equations $MM'=I$ and $M'M=I$.
Schofield~\cite[ch4]{Sch85} gave a slightly more general construction,
inverting a set $\Sigma$ of homomorphisms between finitely
 generated projective $A$-modules.  

Given any ring homomorphism $A\to B$ we may define $\Sigma$ to be the
set of matrices in $A$ which are
invertible in $B$ obtaining 
\begin{equation*}
A\mapright{i_\Sigma} \Sigma^{-1}A \mapright{\gamma} B~.
\end{equation*}
Every matrix with entries in $\Sigma^{-1}A$ can be expressed
(non-uniquely) in the form $f\sigma^{-1}g$ where 
$f$, $\sigma$ and $g$ are matrices with entries in $A$ and
$\sigma\in\Sigma$ (see for example~Schofield~\cite[p52]{Sch85}).

We shall also need the following lemma in
section~\ref{section:Many_Indeterminates}: 
\begin{lemma}
\label{unit_lemma}
A matrix $\alpha$ with entries in $\Sigma^{-1}A$ is invertible
if and only if its
image $\gamma(\alpha)$ is invertible. In particular, 
$\Im(\gamma)^\bullet = B^\bullet \cap \Im(\gamma)$.
\end{lemma} 
\begin{proof}
The `only if' part is easy.
Conversely, suppose $\gamma(\alpha)$ is invertible and
$\alpha=f\sigma^{-1}g$ as above.
The equation
\begin{equation}
\label{eqn:stabilize}
\left(\begin{matrix}
1 & 0 \\
0 & \sigma
\end{matrix}\right)
\left(\begin{matrix}
1 & f \\
0 & 1
\end{matrix}\right)
\left(\begin{matrix}
f\sigma^{-1}g & 0 \\
0 & 1
\end{matrix}\right)
\left(\begin{matrix}
1 & 0 \\
-\sigma^{-1}g & 1
\end{matrix}\right)
=
\left(\begin{matrix}
0 & f \\
-g & \sigma
\end{matrix}\right)
\end{equation}
implies that $\alpha$ is invertible if and only if $\left(\begin{matrix}
0 & f \\
-g & \sigma
\end{matrix}\right)$ is invertible.
But applying $\gamma$ to equation~(\ref{eqn:stabilize}) we learn that 
$\gamma\left(\begin{matrix}
0 & f \\
-g & \sigma
\end{matrix}\right)$ is invertible and hence that 
$\left(\begin{matrix}
0 & f \\
-g & \sigma
\end{matrix}\right)\in \Sigma$.
Thus $\left(\begin{matrix}
0 & f \\
-g & \sigma
\end{matrix}\right)$ and $\alpha$ are invertible over $\Sigma^{-1}A$.
\end{proof}  
\subsection{Proof of Proposition~\ref{Ring_Level_Prop}}
We recall that $S$ denotes the ring
$\Z\langle f,s,g\ |\ fg, fsg, fs^2g, \cdots\rangle$
 and let $\Sigma$ be the set of matrices $\sigma=\sigma_0 + \sigma_1 x
+ \cdots + \sigma_n x^n$ with entries in $S[x]$ such that $\sigma_0$
is invertible (so $\sigma$ is invertible in $S[[x]]$).

We will prove the following two statements:
\begin{itemize}
\item[\ref{Ring_Level_Prop}a)]
The element $f(1-sx)^{-1}g$ is non-zero in $\Sigma^{-1}S[x]$. 
\item[\ref{Ring_Level_Prop}b)]
$f(1-sx)^{-1}g$ lies in the kernel of the natural map
$\gamma:\Sigma^{-1}S[x] \rightarrow S[[x]]$.
\end{itemize}
The second statement \ref{Ring_Level_Prop}b) follows directly from the definition of $S$
\begin{equation*}
\gamma(f(1-sx)^{-1}g)~=~fg + (fsg)x + (fs^2g)x^2 + \cdots~=~0 \in
S[[x]]. 
\end{equation*}
To prove \ref{Ring_Level_Prop}a) we express $S$ once again as the
direct limit $\varinjlim S_m$ with
\begin{equation*}
S_m:=\Z\langle f,s,g\ |\ fg,fsg,\cdots,fs^mg\rangle
\end{equation*}
and the augmentations $\epsilon:S_m[x]\to S_m; x\mapsto 0$ fit into a
commutative diagram 
\begin{equation*}
\xymatrix{
\cdots \ar[r] & S_m[x] \ar[d]^-{\displaystyle{\epsilon}} \ar[r]^-{\displaystyle{p_m}} & S_{m+1}[x] \ar[r]
\ar[d]^-{\displaystyle{\epsilon}}  & \cdots \\
\cdots \ar[r] & S_m \ar[r]^-{\displaystyle{p_m}} & S_{m+1} \ar[r] & \cdots
}  
\end{equation*}
Let $\Sigma_m$ denote the set of matrices
in $S_m[x]$ which become invertible under $\epsilon$, so 
that $p_m(\Sigma_m)\subset \Sigma_{m+1}$ and
$\Sigma = \varinjlim \Sigma_m$.
By lemma~\ref{lemma:localization_directlimits} of appendix~\ref{appendix:directlimits}
\begin{equation*}
\Sigma^{-1}S[x]=\varinjlim \Sigma_m^{-1}S_m[x]
\end{equation*}
so it suffices to show that
$f(1-sx)^{-1}g\neq0 \in \Sigma_m^{-1}S_m[x]$ for each $m\in \N$. 
But $\gamma(f(1-sx)^{-1}g)=\sum^\infty_{n=0} (fs^ng)x^n$ which is
non-zero in $S_m[[x]]$ 
because there does not exist an equation
\begin{equation*}
fs^ng=r_0 fg r'_0 + r_1 fsg r'_1 + \cdots + r_m fs^m g r'_m \in \Z\langle
f,s,g \rangle
\end{equation*}
with $n>m$ and $r_i,r'_i \in \Z\langle f,s,g\rangle$ for $1\leq i \leq
m$. Thus $f(1-sx)^{-1}g\neq 0\in \Sigma_m^{-1}S_m[x]$ and 
the proof of proposition~\ref{Ring_Level_Prop} is complete.

%%%
%%%
\section{Many Indeterminates}
\label{section:Many_Indeterminates}
Let $A$ be any ring, let $X=\{x_1,\cdots, x_\mu\}$ be a finite set,
and let $X^*$ be the free monoid of words in the alphabet $X$. The
free $A$-algebra 
\begin{equation*}
A\langle X\rangle=A\otimes_\Z \Z\langle X\rangle
\end{equation*}
is graded by word length in $X^*$
and is therefore a subring of its completion $A\langle\langle
X\rangle\rangle$ the elements of which are formal power series
$p=\sum_w p_ww$ with $p_w\in A$ for each $w\in X^*$.

Let $\Sigma$ denote the set of matrices in $A\langle X\rangle$ which
are sent to an invertible matrix by the augmentation
$\epsilon:A\langle X\rangle \to 
A;x_i\mapsto 0$ for all $i$. $\Sigma$ is precisely the set
of matrices which are invertible over $A\langle\langle
X\rangle\rangle$ so the inclusion of $A\langle X\rangle$ in
$A\langle\langle X\rangle\rangle$ factors uniquely through $\Sigma^{-1}A\langle
X\rangle$:
\begin{equation*}
A\langle X\rangle \mapright{i_\Sigma} \Sigma^{-1}A\langle X\rangle \mapright{\gamma} A\langle\langle
X\rangle\rangle~.
\end{equation*}
\subsection{Rational Power Series}
\label{subsection:rational_power_series}
In this section we describe the image of $\gamma$.
\begin{definition}
Let $\CR^A$ denote the rational closure of $A\langle X\rangle$. 
In other words $\CR^A$ is the intersection of all the rings $R$ such that
$A\langle X\rangle \subset R\subset A\langle\langle X\rangle\rangle$
and $R^\bullet=R\cap A\langle\langle X\rangle\rangle^\bullet$.
A power series $p\in \CR^A$ is said to be {\it rational}.
\end{definition}
\begin{proposition}
\label{image_of_gamma}
$\gamma(\Sigma^{-1}A\langle X\rangle)=\CR^A.$
\end{proposition}
\begin{proof}
To prove $\CR^A\subset \Im(\gamma)$, we need only note
that $\Im(\gamma)^\bullet = \Im(\gamma)\cap A\langle\langle
X\rangle\rangle^\bullet$ by lemma~\ref{unit_lemma}
above. 

Conversely, to prove that $\Im(\gamma)\subset \CR^A$ 
it suffices to show that every matrix
$\sigma\in\Sigma$ has an inverse with entries in $\CR^A$
so that there is a commutative diagram
\begin{equation*}
\xymatrix{
A\langle X\rangle \ar[d]_-{\displaystyle{i_\Sigma}} \ar[dr] \ar[drr] & & \\
\Sigma^{-1}A\langle X\rangle \ar[r] & \CR^A \ar[r] & A\langle\langle
X\rangle\rangle.
}
\end{equation*} 
Recall that $\epsilon:A\langle X\rangle \to A$ is the augmentation
given by $\epsilon(x_i)=0$ for all $i$.
Multiplying $\sigma$ by $\epsilon(\sigma)^{-1}$ if necessary we can
assume that $\epsilon(\sigma)=I$.
Each diagonal entry of $\sigma$ has an inverse in $\CR^A$ so, after
elementary row operations (which are of course invertible), $\sigma$ becomes
a diagonal matrix where each diagonal entry~$\sigma_{ii}$ has
$\epsilon(\sigma_{ii})=1$ (cf~the proof of lemma~\ref{lemma:K1powerseries}i)
above). 
By the definition of~$\CR^A$, each 
$\sigma_{ii}$ has an inverse in~$\CR^A$.
\end{proof}
\subsection{Sch\"utzenberger's Theorem}
\begin{proposition}
\label{linear_normal_form}
Every matrix $\alpha$ with entries in $\Sigma^{-1}A\langle X\rangle$
can be expressed (non-uniquely) in the form 
\begin{equation}
\label{equation:linear_normal_form}
\alpha=f(1-s_1x_1-\cdots-s_\mu x_\mu)^{-1}g
\end{equation}
where $f$, $s_1$,\ldots, $s_\mu$ and  $g$ are matrices with entries in~$A$.
\end{proposition}
\begin{proof} 
It suffices to show that $\alpha$ has the form $f\sigma^{-1}g$ where
$f$ and $g$ have entries in~$A$ and $\sigma$ is a linear matrix 
$\sigma=\sigma_0+\sum_{i=1}^\mu \sigma_ix_i$ with $\sigma_0$ invertible.
For then $\alpha=(f\sigma_0)(\sigma_0^{-1}\sigma)g$.

Note first that if $\alpha_1=f_1\sigma_1^{-1}g_1$ and
$\alpha_2=f_2\sigma_2^{-1}g_2$ then 
\begin{align}
\alpha_1 - \alpha_2 &= \left(
\begin{matrix}
f_1 & -f_2
\end{matrix}\right)
\left(\begin{matrix}
\sigma_1 & 0 \\
0 & \sigma_2
\end{matrix}\right)^{-1}
\left(\begin{matrix}
g_1 \\
g_2
\end{matrix}\right)
\label{equation:subtraction}
\\
\intertext{and}
\alpha_1\alpha_2 &=
\left(\begin{matrix}
f_1 & 0
\end{matrix}\right)
\left(\begin{matrix}
\sigma_1 & -g_1f_2 \\
0 & \sigma_2
\end{matrix}\right)^{-1}
\left(\begin{matrix}
0 \\
g_2
\end{matrix}\right)
\label{equation:composite}
\end{align}
whenever the left-hand sides make sense (cf~\cite[p52]{Sch85}).
Hence, we need only treat the cases where
i) $\alpha$ has entries in $A\langle X\rangle$ and
ii) $\alpha=\sigma^{-1}$ with  $\sigma\in \Sigma$.

If $\alpha$ has entries in $A\langle X\rangle$ then by repeated application
of the equation 
\begin{equation}
\label{equation:linearising_step}
\left(\begin{matrix}
a+bc & 0 \\
0 & 1
\end{matrix}\right)
=
\left(\begin{matrix}
1 & -b \\
0 & 1
\end{matrix}\right)
\left(\begin{matrix}
a & b \\
-c & 1
\end{matrix}\right)
\left(\begin{matrix}
1 & 0 \\
c & 1
\end{matrix}\right),
\end{equation}
in which $a$, $b$, $c$ and $1$ denote matrices, 
some stabilisation $\left(\begin{matrix} \alpha & 0 \\
0 & 1 \end{matrix}\right)$ can be expressed as a product of linear matrices. 
Each linear matrix $a_0+a_1x_1+\cdots + a_\mu x_\mu$ can be written
\begin{equation*}
\left(\begin{matrix}
1 & 0
\end{matrix}\right)
\left(\begin{matrix}
1 & -a_0 \\
0 & 1
\end{matrix}\right)^{-1}
\left(\begin{matrix}
0 \\ 
1
\end{matrix}\right)
~+~
\sum_{i=1}^\mu
\left(\begin{matrix}
1 & 0
\end{matrix}\right)
\left(\begin{matrix}
1 & -a_ix_i \\
0 & 1
\end{matrix}\right)^{-1}
\left(\begin{matrix}
0 \\
1
\end{matrix}\right)
\end{equation*}
and equations~(\ref{equation:subtraction}) and~(\ref{equation:composite}) imply that 
$\alpha=\left(\begin{matrix}1 &
0\end{matrix}\right)\left(\begin{matrix} \alpha & 0 \\
0 & 1 \end{matrix}\right)\left(\begin{matrix} 1 \\ 0
\end{matrix}\right)$ 
is of the required form~$f\sigma^{-1}g$.

The case $\alpha=\sigma^{-1}$ is similar but slightly easier; we
repeatedly apply equation~(\ref{equation:linearising_step}) to express
(a stabilisation of) $\sigma^{-1}$ as
a product of inverses of linear matrices in $\Sigma$ and then
apply equation~(\ref{equation:composite}). 
\end{proof}

A power series $p\in A\langle\langle X\rangle\rangle$ is said to be
{\it recognisable} if it is of the form
\begin{equation*}
p=fg+\sum_{i=1}^\mu fs_igx_i +
\sum_{i,j=1}^\mu fs_is_jgx_ix_j + \cdots~.
\end{equation*}
where $f\in A^n$ is a row vector, $g\in A^n$ is a column vector and
each $s_i$ is an $n\times n$ matrix in $A$.
Propositions~\ref{linear_normal_form} and~\ref{image_of_gamma} imply
\begin{corollary}[Sch\"utzenberger's theorem]
\label{theorem:rat=rec}
A power series $p\in A\langle\langle X\rangle\rangle$ is rational
if and only if it is recognisable.
\end{corollary}
\section{Localization of the Free Group Ring}
\label{section:Farber_Vogel}
We identify the localization of the group ring of the
free group studied by Farber and Vogel~\cite{FarVog92} with the localization
$\Sigma^{-1}A\langle X\rangle$ of the present paper. 

Let $F_\mu$ denote the
free group on generators $z_1,\dots, z_\mu$ and as usual let $A$
be a (not necessarily commutative) ring.
$AF_\mu$ will denote the group ring, in which the elements of the
group $F_\mu$ are assumed to commute with elements of $A$.
Let $\epsilon:AF_\mu\rightarrow A;\ z_i\mapsto 1$ for all $i$ and
let $\Psi$ denote the set of square matrices $M$ in $AF_\mu$
such that $\epsilon(M)$ is invertible. $\Psi$ is
denoted $\Sigma$ in~\cite{FarVog92}.

All the matrices in $\Psi$ become invertible under the Magnus embedding of
the group ring
\begin{align*}
AF_\mu &\rightarrow A\langle\langle X\rangle\rangle \\
z_i &\mapsto 1 + x_i \\
z_i^{-1} &\mapsto 1-x_i + x_i^2 - x_i^3 \cdots.
\end{align*}
so the embedding factors through $\Psi^{-1}AF_\mu$
\begin{equation*}
AF_\mu\mapright{i_\Psi} \Psi^{-1}AF_\mu \mapright{\gamma}
A\langle\langle X\rangle\rangle~.
\end{equation*}
Farber and Vogel proved that if $A$ is a (commutative) principle ideal
domain then $\gamma$ is an injection and the image of $\gamma$ is the
ring $\CR^A$ of rational power series.

For any ring~$A$ let $m:A\langle X\rangle \to AF_\mu$ be the ring homomorphism
defined by $x_i\mapsto z_i-1$ for all $i$. There is
a commutative diagram
\begin{equation*}
\xymatrix{
A\langle X\rangle \ar[r]^-{\displaystyle{\epsilon}}
\ar[d]_-{\displaystyle{m}} & A \\
AF_\mu \ar[ur]_-{\displaystyle{\epsilon}} & \\
}
\end{equation*}
so $m(\Sigma)\subset \Phi$ 
and $m$  
induces a homomorphism $m:\Sigma^{-1}A\langle X\rangle \to
\Psi^{-1}AF_\mu$ which fits into a commutative diagram
\begin{equation*}
%\label{diagram:groupring}
\xymatrix{
A\langle X\rangle \ar[r]^{i_\Sigma} \ar[d]_m &
 \Sigma^{-1}A\langle X\rangle \ar[d]^m \ar[r]^\gamma &  A\langle\langle
X\rangle\rangle \\ 
AF_\mu \ar@{-->}[ur]^l \ar[r]^{i_\Psi} & \Psi^{-1}AF_\mu \ar[ur]_\gamma & }
\end{equation*}
\begin{proposition}
$m:\Sigma^{-1}A\langle X\rangle \to
\Psi^{-1}AF_\mu$ is an isomorphism.
\end{proposition}
\begin{proof}
Observe first that $AF_\mu$ is isomorphic to the Cohn localization 
\begin{equation*}
\{1+x_i\mid 1\leq i \leq\mu\}^{-1}A\langle X\rangle
\end{equation*}
inverting the $1\times1$ matrices $(1+x_i)$. 
Since $(1+x_i)\in\Sigma$ the homomorphism~$i_\Sigma$ factors
uniquely through $AF_\mu$ as indicated by the broken arrow $l$ in the commutative
diagram above. Explicitly, $l:AF_\mu\to \Sigma^{-1}A\langle X\rangle;
z_i\mapsto 1+x_i$. Now if $\psi\in\Psi$ then $\gamma l(\psi)$ is invertible so 
by lemma~\ref{unit_lemma} $l(\psi)$ is invertible. Thus $l$ induces a
map $\Psi^{-1}AF_\mu\to \Sigma^{-1}A\langle X\rangle$ which, by the
universal properties of $i_\Sigma$ and $i_\Psi$, 
is inverse to $m$.
\end{proof}
%%%
%%%
\appendix
\section{Direct Limits}
\label{appendix:directlimits}
In this appendix we prove that Cohn localization and the functor
$\rEnd(\functor)$ commute with direct limits. 
\subsection{Cohn Localization}
First we make the former claim more precise. 
Suppose $I$ is a directed set and $(\{A_m\}_{m\in I},\{f_m^{l}:A_m
\to A_{l}\}_{m\leq l})$ is a direct system of rings.
Suppose further that for each $m\in I$ we have a set of matrices
$\Sigma_m$ with entries in $A_m$ such that $f_m^{l}(\Sigma_m)\subset
\Sigma_{l}$ whenever $m\leq l$.   
If $i_m:A_m\to \Sigma_m^{-1}A_m$ is the
universal $\Sigma_m$-inverting ring homomorphism for each $m$, then
when $m\leq l$ the composite
\begin{equation*}
A_m\mapright{f_m^{l}}A_{l}\mapright{i_{l}} \Sigma_{l}^{-1}A_{l}
\end{equation*}
is $\Sigma_m$-inverting and therefore factors through a map 
$\Sigma^{-1}f_m^{l}:\Sigma_m^{-1}A_m\to \Sigma_{l}^{-1}A_{l}$.
It is easy to see that
$\Sigma^{-1}f_{l}^{k}\circ \Sigma^{-1}f_m^{l}=\Sigma^{-1}f_m^{k}$ when
$m\leq l \leq k$.

For any ring $A$ let $\CM(A)$ denote the set of matrices (of any size
and shape) with entries in $A$. 
The inclusions $\Sigma_m\subset \CM(A_m)$
induce an injection
\begin{equation*}
\varinjlim \Sigma_m \to \varinjlim \CM(A_m) = \CM(\varinjlim A_m).
\end{equation*}
\begin{lemma}
\label{lemma:localization_directlimits}
  There is a natural isomorphism
\begin{equation*}
(\varinjlim \Sigma_m)^{-1} (\varinjlim A_m) \cong \varinjlim
(\Sigma_m^{-1}A_m).
\end{equation*}
\end{lemma} 
\begin{proof}
One can check that the canonical map 
%\begin{equation*}
$\varinjlim i_m : \varinjlim A_m \to \varinjlim \Sigma_m^{-1}A_m$
%\end{equation*}
 is universal among $(\varinjlim \Sigma_m)^{-1}$-inverting
 homomorphisms.
The details are left to the reader.
\end{proof}
\subsection{The Endomorphism Class Group}
\label{subappendix:end_dir_lim}
%We consider next the functor $\rEnd_0(\functor)$. 
\begin{lemma}
\label{lemma:rEnd_directlimits}
There is a natural isomorphism
\begin{equation*}
\varinjlim \rEnd_0(A_m) \cong \rEnd_0(\varinjlim A_m).
\end{equation*}
\end{lemma}
\begin{proof}
The canonical maps $f_m:A_m\to \varinjlim A_m$ induce maps 
$f_m:\rEnd_0(A_m)\to \rEnd_0(\varinjlim A_m)$ satisfying
$f_{l}f_m^{l} = f_m$ for $m\leq l$. We aim to
prove that any other system of maps $g_m :
\rEnd_0(A_m)\to T$ with $g_{l}f_m^{l} = g_m$ for $m\leq l$ factors
uniquely through $\rEnd_0(\varinjlim A_m)$:
\begin{equation*}
\xymatrix{
\{\rEnd_0(A_m)\} \ar[rr]^{\{g_m\}} \ar[d]^{\{f_m\}} & & T \\
\rEnd_0(\varinjlim A_m) \ar@{-->}[urr]^{g} &
}
\end{equation*}

Suppose $[M]$ is a generator of $\rEnd_0(\varinjlim A_m)$ where $M\in
M_n(\varinjlim A_m)$. $M$ is the image $f_m(M_m)$ of some matrix $M_m\in
M_n(A_m)$ so we can define $g[M]=g_m[M_m]$. To show $g$
is well-defined there are two things to check: \smallskip \\ \noindent
i) If $M_{l}\in M_n(A_l)$ is an alternative choice with $f_l(M_l)=M$
then we require $g_m[M_m]=g_l[M_l]$. Indeed, there exists $k$ such
that $l\leq k$, $m\leq k$ and $f_l^k(M_k)=f_m^k(M_m)\in M_n(A_k)$. Hence
$g_m[M_m]=g_kf_m^k[M_m]=g_kf_l^k[M_l]=g_l[M_l]$. \smallskip \\ \noindent
ii) We must check that $g$ respects the defining relations of
$\rEnd_0(\varinjlim A_m)$. 
\begin{enumerate}
\item A matrix $\left(\begin{matrix}
M & N \\
0 & M'
\end{matrix}\right)$
is the image of some matrix $\left(\begin{matrix}
M_m & N_m \\
0 & M'_m
\end{matrix}\right)$
so \\ \noindent
$g\left[\begin{matrix}
M & N \\
0 & M'
\end{matrix}\right] = g_m\left[\begin{matrix}
M_m & N_m \\
0 & M'_m
\end{matrix}\right]=g_m([M_m]+[M'_m])=g[M]+g[M']$.
\item Suppose $M'=PMP^{-1}$ for some invertible matrix $P$.
For large enough $m$ we can choose $P_m,Q_m\in M_n(A_m)$ to represent
$P$ and $P^{-1}$ respectively. Since
$I=f_m(P_m)f_m(Q_m)$ there exists $k\geq m$ such that 
$P_kQ_k=I\in M_n(A_k)$ where $P_k=f_m^kP_m$ and $Q_k=
f_m^kQ_m$. Thus $M'=f_k(P_kM_kP_k^{-1})$ and
$g[M']=g_k[P_kM_kP_k^{-1}]=g_k[M_k]=g[M]$.
\item If $M$ is the zero matrix, $g[M]=0$.
\end{enumerate}
Uniqueness of $g$ follows from the fact that every class $[M]$
in $\rEnd_0(\varinjlim A_m)$ is an image of a class $[M_m]\in\rEnd_0(A_m)$.
\end{proof}

\bigskip

\noindent
Department of Mathematics and Statistics, \\
The University of Edinburgh, \\
King's Buildings, \\
Edinburgh EH9 3JZ, \\
Scotland, UK. \bigskip \\
e-mail: des@sheiham.com
\end{document}